\documentclass[11pt]{article}
\usepackage[margin=1in]{geometry}

\usepackage{setspace}
\doublespace

\usepackage{amsmath}
\usepackage{amssymb}
\usepackage{bm}
\newtheorem{theorem}{Theorem}[section]
\newtheorem{proposition}{Proposition}[section]
\newtheorem{corollary}{Corollary}[section]
\newtheorem{lemma}{Lemma}[section]

\newcommand{\RR}{\mathbb{R}}

\newcommand{\ZZ}{\mathbb{Z}}

\newcommand{\Pois}{\mathrm{Pois}}
\newcommand{\eL}{{\cal L}}
\newcommand{\exponent}[1]{\exp\{#1\}}           % exponent
\newcommand{\Exponent}[1]{\exp\Bigl\{#1\Bigr\}} % large exponent
\newcommand{\norm}[1]{\|#1\|}               % norm
\newcommand{\ab}[1]{\vert#1\vert}           % absolut value
\newcommand{\Ab}[1]{\Big\vert#1\Big\vert}   % large absolute value
\newcommand{\qed}{$\square$}                % end of proof
\newcommand{\Expect}{\mathrm{E}}            % expectation
\newcommand{\EE}{\mathrm{E}}            % expectation
\newcommand{\ii}{{\mathrm i}}
\newcommand{\ee}{{\mathrm e}}
\newcommand{\dd}{{\mathrm{d}}}

\newcommand{\boldt}{\mathbf{t}}
\newcommand{\boldzero}{\mathbf{0}}
\newcommand{\boldm}{\mathbf{m}}

\newcommand{\bolda}{\mathbf{a}}
\newcommand{\boldp}{\mathbf{p}}
\newcommand{\boldbeta}{\bm{\beta}}

\newcommand{\w}{\widehat}
\newcommand{\vfi}{\varphi}
\newcommand{\wE}{\w\Expect}
\newcommand{\veps}{\varepsilon}

 \scrollmode

\begin{document}
%%%%%%%%%%%%%%%%%%%%%%%%%%%%%%%%%%%%%%%%%%%%%%%%%%%%%%%%%%%%%%%%%%%%%%
%%%%%%%%%%%%%%%%%%%%%%%%%%%%%%%%%%%%%%%%%%%%%%%%%%%%%%%%%%%%%%%%%%%%%%
\title{
%%\textbf{\small\fbox{
%%\begin{minipage}{12cm}
%%\begin{center}
%%Filename/Date(day.month.year)/Time of
%%compilation of dvi-file: \\
%%\jobname.tex/\mydate\today/\currenttime
%%\end{center}
%%\end{minipage}
%%}}\vspace{5mm}\\
%Multivariate Poisson type approximations\\
 Compound Poisson approximations in $\ell_p$-norm for sums of weakly dependent  vectors  }
\author{V. \v{C}ekanavi\v{c}ius \footnote{
Corresponding author, postal address : Department of Mathematics
and Informatics, Vilnius University, Naugarduko 24, Vilnius
LT-03225, Lithuania. E-mail: vydas.cekanavicius@mif.vu.lt} \;\;and\, P. Vellaisamy
\footnote{ Postal address:  Department of Mathematics, Indian Institute of Technology Bombay, Powai, Mumbai-
400076, India. E-mail: pv@math.iitb.ac.in}
\\
}
\date{\today}
\maketitle
%%%%%%%%%%%%%%%%%%%%%%%%%%%%%%%%%%%%%%%%%%%%%%%%%%%%%%%%%%%%%%%%%%%%%%
%%%%%%%%%%%%%%%%%%%%%%%%%%%%%%%%%%%%%%%%%%%%%%%%%%%%%%%%%%%%%%%%%%%%%%
\begin{abstract}

The distribution of the sum of 1-dependent lattice vectors with supports
on coordinate axes is approximated by a multivariate compound
Poisson distribution and by signed compound Poisson measure. The  local  and $\ell_\alpha$-norms are used to obtain the error bounds.
The Heinrich method is used for the proofs.

{\emph{Keywords: }\small    % arranged in alphabetical order
Compound  Poisson distribution;  expansion in the exponent; $\ell_p$ norm; local norm;
 multivariate  distribution;
\medskip\\
} {\small MSC 2000 Subject Classification:
primary  62E17; %Approximations to distributions (nonasymptotic)
secondary 60F05.   % Central limit and other weak theorems
}
\end{abstract}
%\section{Introduction}
\label{intro}

Numerous papers are devoted to Poisson and compound Poisson approximations in one-dimensional case, see, for example,  surveys \cite{ChR13,Nov19}.
 Various metrics, such as local, Gini (Wasserstein) metric, chi-square metric and analogues of $\ell_p$ norms were used, see, for example, \cite{BGX15,BRR19,HZ10,N18,ZH10}.

The  multivariate case is less explored.
 For compound Poisson approximations in L\'evy, L\'evy - Prokhorov and Kolmogorov metrics, see \cite{zc89,zc92,zc14}.
 Multivariate Poisson approximation in total variation for sums of independent lattice
 vectors  concentrated on  coordinate  vectors is considered in
  \cite{Ar76,barb05,roo99}, compound Poisson and signed compound Poisson approximations are applied in \cite{KC14,roo03,Roos04,roo17}.
 In \cite{KC14,roo99}  local point metric is also used. Necessary and sufficient conditions for the weak convergence of the distribution of a sum of dependent random vectors to a compound Poisson vector are given in \cite{N11}, Theorems 6.3 and 6.4.  An estimate  of Compound Poisson approximation to the sum of weakly dependent discrete random vectors is given in \cite{N11}, Theorem 6.8.

In this paper, we similarly investigate the sum of random vectors (rvs) concentrated on $k$-dimensional unit vectors. We assume that rvs are 1-dependent and  estimate the accuracy of  approximations in  local and $\ell_p$ norms.
As usual 1-dependence means that sigma-algebras generated by $X_1,X_2,\dots,X_k$ and $X_t, X_{t+1},\dots, X_s$  are independent for any $k-t>1$.
This implies for example that $X_1$ and $X_3$, $X_1$ and $X_4$ and etc.  are independent.  Note that any sum of m-dependent dependent rvs can be reduced to the sum of 1-dependent rvs by grouping of consequent summands.

Let $\boldzero=(0,\dots,0)$ and let $e_r$ denote the $r$'th coordinate vector in $\RR^k$, that is,
$e_r=(0,\dots,1,\dots,0)$, $ 1 \leq r \leq k.$
Further, let $X_1, X_2,\dots, X_n$ be 1-dependent identically distributed
 $k$-dimensional rvs and
$P(X_1=e_r)=p_r$, $p_r\in (0,1)$, $r=1,2,\dots,k$, $P(X_1=\boldzero)=1-(p_1+p_2+\cdots+p_k)$.
The dependence of consequent summands is reflected in joint probabilities $p_{rj}=P(X_1=e_r, X_2=e_j)$,
$p_{rjm}=P(X_1=e_r,X_2=e_j,X_3=e_m)$.  The distribution of $S_n=X_1+X_2+\cdots+X_n$ is denoted by $F_n$.

Our aim is to approximate the distribution of $S_n$.
Approximations used in this paper are mostly defined by their Fourier transforms.
If a measure
$M$ is concentrated on $k$-dimensional set of integers $\ZZ^k=\ZZ\times\ZZ\times\cdots\times\ZZ$, then its Fourier transform
(characteristic function)  is denoted by
\begin{equation}\label{Ftrans}
\w M(\boldt)=\sum_{\boldm\in\ZZ^k}M\{\boldm\}\exponent{\ii(\boldm, \boldt)}.\end{equation}
Here and henceforth $(\boldm,\boldt):=m_1t_1+m_2t_2+\cdots+m_kt_k$ and $\ii$ is imaginary unit.
Note that there is one to one correspondence between a distribution and its characteristic function. Observe also  that $\w F_n(\boldt)=\EE\exponent{\ii (S_n,\boldt)}$ and $\EE\ee^{\ii(X_1,\boldt)}=1+\sum_{r=1}^kp_r(\ee^{\ii t_r}-1)$.

%In general, k-dimensional Poisson distribution is defined as a distribution %of  random sum $Y_1+Y_2+\cdots+Y_N$, where $N$ is
 %Poisson random variable  with parameter $n$, and $Y_1,Y_2,\dots$ are %independent identically distributed rvs concentrated
% on vectors with components 0 or 1 (such as $(1,1,0,...)$ etc.) It is also %assumed that $Y_1,Y_2,\dots$ do not depend on $N$.

We say rv $\tilde Y=(\tilde Y_1, \tilde Y_2, \cdots, \tilde Y_k)$ follows a $k$-dimensional Poisson
if $\tilde Y_j$'s are independent and  $\tilde Y_j$ follows the Poisson distribution with
parameter $\mu_j$, $ 1 \leq j \leq k$. It is denoted by $\Pois(\bm{\mu})$, $\bm{\mu}=(\mu_1, \mu_2,\dots, \mu_k)$.  In this paper, we consider the problem of approximating
the distribution of $S_n$ to the $ k$-dimensional Poisson distribution  $\Pois(\bm{\lambda})$, where $\bm{\lambda}=(np_1,np_2,\dots,np_k)$. Let $P(\lambda,s)=\ee^{=\lambda}\lambda^s/s!$ denote probability of Poisson random variable, Then
\[\Pois(\bm{\lambda})\{\boldm\}=\prod_{j=1}^kP(np_j,m_j),\quad\w\Pois(\bm{\lambda})(\boldt)=\Exponent{n\sum_{j=1}^kp_j(\ee^{\ii t_j}-1)}.\]

First order asymptotic expansion is constructed as $\Pois(\bm{\lambda})+A_1$, where measure $A_1$ has the following Fourier transform:
\begin{eqnarray}
\w A_1(\boldt)&:=&\w\Pois(\bm{\lambda})(\boldt)\bigg(-\frac{n}{2}\bigg(\sum_{j=1}^kp_j(\ee^{\ii t_j}-1)\bigg)^2\nonumber\\
&&+
(n-1)\sum_{j,m=1}^k(p_{jm}-p_jp_m)(\ee^{\ii t_j}-1)(\ee^{\ii t_m}-1)\bigg).
\label{APois}
\end{eqnarray}
Let $\Delta P(\lambda,s):=P(\lambda,s)-P(\lambda,s-1)$, $\Delta^2P(\lambda,s)=\Delta(\Delta P(\lambda,s))$ and let $d_{jr}=(n-1)(p_{jr}-p_jp_r)-(n/2)p_jp_r$. From the formula of inversion it follows that, for any $\boldm\in\ZZ^k$,
\begin{eqnarray*}
A_1\{\boldm\}&=&\sum_{j=1}^kd_{jj}\Delta^2 P(np_j,m_j)\prod_{l\ne j}^kP(np_l,m_l)
\\
&&+\sum_{j,r=1; j\ne r}^kd_{jr}\Delta P(np_j,m_j)\Delta P(np_r,m_r)\prod_{l\ne j,r}P(np_l,m_l).
\end{eqnarray*}

If we construct similar asymptotic expansion in the exponent, then the result is signed compound Poisson measure $G$ with Fourier transform
\begin{eqnarray}
\w G(\boldt)&:=&\w\Pois(\bm{\lambda})(\boldt)\exp\bigg\{-\frac{n}{2}\bigg(\sum_{j=1}^kp_j(\ee^{\ii t_j}-1)\bigg)^2\nonumber\\
&&+
(n-1)\sum_{j,m=1}^k(p_{jm}-p_jp_m)(\ee^{\ii t_j}-1)(\ee^{\ii t_m}-1)\bigg\}.\label{wG}
\end{eqnarray}
Formula of inversion also allows  explicit expression for $G\{\boldm\}$. However, the resulting formula is quite long and plays no role in further proofs, therefore, is omitted. %left as an exercise for the reader.
The idea to use compound Poisson type signed measures by retaining a part of asymptotic expansion in the exponent goes back  to early eighties of the twentieth century. Inspite of obvious structural similarity, generally such measures ensure much better accuracy than Poisson or even the second-order Poisson approximations,
see for example, \cite{CeVe15,Nov19,roo03} and the references therein. Both the measures $A_1$ and $G$ can be written as convolutions of measures concentrated on various $e_r$.

In this paper, symbol $p$ is reserved for probabilities. Therefore, in the definition of norms, we instead use symbol $\alpha$, which is further on  assumed to be fixed positive number.
We define respectively local and $\ell_\alpha$-norms for finite measure $M$  concentrated on $\ZZ^k$ as
\[
\quad \norm{M}_\infty=\sup\limits_{\boldm\in\ZZ^k}\ab{M\{\boldm\}},\quad \norm{M}_\alpha=\bigg(\sum_{\boldm\in\ZZ^k}M^\alpha\{\boldm\}\bigg)^{1/\alpha},\quad \alpha\in [1,\infty).\]
The case $\alpha=1$ corresponds to the total variation norm $\norm{M}:=\norm{M}_1$. Local norm can be viewed as a limit case of $\norm{M}_\alpha$ when $\alpha\to\infty$.
Thus, total variation and local norms form natural boundaries for all $\ell_\alpha$  norms. In this paper, the emphasis is on local and $\ell_\alpha$, ($\alpha\geq 2$) norms. Note that total variation norm is equivalent  to the total variation distance. More precisely, $d_{TV}(F,G):=\sup\ab{F\{B\}-G\{B\}}=\frac{1}{2}\norm{F-G}$. Here supremum is taken over all $k$-dimensional Borel sets $B$.

 We
denote by $C$ positive absolute constants, the values of which may
change from line to line, or even within the same line. Similarly,
by $C(\cdot)$ we denote constants depending on the indicated
argument only. Sometimes, to avoid possible ambiguity we supply
$C$ with index.
 Similarly, $\theta$ is used for a real or a complex number satisfying $\ab{\theta}\leq 1$.

%%%%%%%%%%%%%%%%%%%%%%%%%%%%%%%%%%%%%%%%%%%%%%%%%%%%%%%%%%%%%%%%%%%%%%
%%%%%%%%%%%%%%%%%%%%%%%%%%%%%%%%%%%%%%%

%%%%%%%%%%%%%%%%%%%%%%%%%%%%%%%%%%%%%%%%%%%%%%%%%%%%%%%%%%%%%%%%%%%%%%%%%%%%%%%%%%%%%%%%%%%%%%%%%%%%%%%%%%%%%%%%%%%%%%%%
%%%%%%%%%%%%%%%%%%%%%%%%%%%%%%%%%%%%%%%%%%%%%%%%%%%%%%%%%%%%%%%%%%%%%%%%%%%%%%%%%%%%%%%%%%%%%%%%%%%%%%%%%%%%%%%%%%%%%%%%
\section{Some known results}

The  most part of multivariate results related to compound Poisson  approximations are proved for independent rvs concentrated on
 $e_r, 1 \leq r \leq k.$  Set $p=p_1+p_2+\cdots+p_k$.
The total variation metric  is typically used to estimate the accuracy of approximations. If we assume that $\tilde X_j$ are independent copies of $X_j$, $j=1,2,\dots,n$,
and denote the distribution of $W_n= \tilde X_1+\tilde X_2+\dots+\tilde X_n$ by
 $\eL(W_n)$,  then
\begin{equation*}
\frac{1}{7}\max_{1\leq r\leq k}\min(p_r,np_r^2)\leq\norm{\eL(W_n)-\Pois(\bm{\lambda})}\leq 2\min(p,np^2).\label{trick}
\end{equation*}
The upper bound  follows directly from one-dimensional Poisson approximation to the binomial distribution;
see \cite{Ce16}, p.~29 and Equation (1.1) in \cite{roo17}. The lower bound  is  a special case of Proposition 1.3 in  \cite{roo17}. Note that constant 2 on the right-hand side of (\ref{trick}) is not the optimal one and there exist other, longer expressions, with smaller constants see, for example Eq. (50)-(51) in \cite{Nov19} or (1.3) in \cite{roo17} or discussion in \cite{N18}.
  In \cite{roo01}  $\ell_\alpha$-norm was used for Krawtchouk expansions, though we are unaware about any similar multivariate Poisson approximation result. In one-dimensional case, squared $\ell_2$-norm was used in the seminal paper of Franken \cite{Fr64} and the closeness of binomial and Poisson distributions was thoroughly investigated in \cite{HZ10} for an analogue of $\ell_\alpha$-norm for even more general case of $\alpha\in (0,\infty)$.

 There are many local  estimates for Poisson approximation to sums of random variables, see for example \cite{BJ89}. However, unlike in (\ref{trick}),  they can not be used for obtaining local estimate for $k$-dimensional vectors. The local  bound follows from Equation (26)  in \cite{roo99}. Let
 \[v(r)=2np_r^2\min\Big\{\frac{1}{2np_r}, \ee\Big\}.\]
If  $\sum_{r=1}^{k} \sqrt{2 v(r)}<1$, then
 \begin{equation}
\norm{\eL(W_n)-\Pois(\bm{\lambda})}_\infty\leq 2\prod_{j=1}^k\min\Big\{\frac{1}{2np_j}, \ee\Big\}\frac{\Big(\sum_{r=1}^k\sqrt{v(r)}\Big)^2}{1-\sum_{r=1}^k\sqrt{2v(r)}}.
 \label{ros1}
 \end{equation}

 If $1-p>C>0$ and all $np_r \geq 1$, then   (\ref{ros1}) implies
 \begin{equation}
 \norm{\eL(W_n)-\Pois(\bm{\lambda})}_\infty\leq C(k)\prod_{j=1}^k\frac{1}{\sqrt{np_j}}\bigg(\sum_{r=1}^k\sqrt{p_r}\bigg)^2.
  \label{ros2}
 \end{equation}

Approximation of sums of independent non-identically distributed  vectors by analogues of $G$ was thoroughly investigated by Roos \cite{Roos04}. For small probabilities a different choice of parameters for $G$ was proposed by Borovkov \cite{B88}. However, arguably the best estimate for identically distributed independent vectors follows from one-dimensional estimate  from \cite{ZH10}: let $\omega:=\sum_{r=1}^kp_r^2/p<1/2$, then
\begin{equation}
	\norm{\eL(W_n)-G}\leq \frac{p^{3/2}}{\sqrt{n}}\bigg(\frac{\sqrt{6}\cdot 0.374}{(1-\omega)^2}+\frac{\sqrt{3\omega}}{2\sqrt{2}(1-\omega)^{5/2}}.\label{ZHwang}\bigg).
\end{equation}
Observe that the accuracy in (\ref{ZHwang}) is at least of the order $O(n^{-1/2})$. Note also that, due to independence, $p_{rm}-p_rp_m=0$ and $G$ has simpler structure than in (\ref{wG}).

From Corollary 6.9 in \cite{N11} it follows that for sums of 1-dependent rvs Poisson approximation in total variation can be of order $O(p\sqrt{n})$, which is significantly weaker than (\ref{trick}) for moderate $p$. To the best of our knowledge, there is no similar result for the 1-dependent vectors in $\ell_\alpha>1$. Therefore, to get the general idea about what can be expected, we formulate two one-dimensional local estimates, which can be easily obtained from Lemmas 6.3, 6.4, 6.7 in \cite{CeVe15} and the inversion formula  (Lemma \ref{locinv}) given in Section 4.
Let $k=1$ and $p_1<0.01$ and $p_{12}<0.05p_1$, $a:=\max(np_1,1)$. Then

 \begin{equation}
\norm{F_n-\Pois(np_1)}_\infty\leq \frac{Cn(p_{12}+p_1^2)}{a\sqrt{a}}
 \label{ve1}
 \end{equation}
and
 \begin{equation}
\norm{F_n-G}_\infty\leq \frac{Cn(p_{123}+p_{12}p+p_1^3)}{a^2}.
  \label{ve2}
 \end{equation}
 Observe  summands $p_{12}$ and $p_{123}$ reflecting the possible 1-dependence of random variables.

%%%%%%%%%%%%%%%%%%%%%%%%%%%%%%%%%%%%%%%%%%%%%%%%%%%%%%%%%%%%%%%%%%%%%%%%%%%%%%%%%%%%%%%%%%%%%%%%%%%%%%%%%%%%%%%%%%%%%%%%

%%%%%%%%%%%%%%%%%%%%%%%%%%%%%%%%%%%%%%%%%%%%%%%%%%%%%%%%%%%%%%%%%%%%%%%%%%%%%%%%%%%%%%%%%%%%%%%%%%%%%%%%%%%%%%%%%%%%%%%%
%%%%%%%%%%%%%%%%%%%%%%%%%%%%%%%%%%%%%%%%%%%%%%%%%%%%%%%%%%%%%%%%%%%%%%%%%%%%%%%%%%%%%%%%%%%%%%%%%%%%%%%%%%%%%%%%%%%%%%%%

\section{Results}

%%%%%%%%%%%%%%%%%%%%%%%%%%%%%%%%%%%%%%%%%%%%%%%%%%%%%%%%%%%%%%%%%%%%%%%%%%%%%%%%%%%%%%%%%%%%%%%%%%%%%%%%5
We begin from assumptions on the smallness of probabilities and their interdependency. The magnitude of constants is determined by the method of proof. We assume that
\begin{equation}
\max_{1\leq j\leq k}p_j\leq \frac{1}{144k},\quad \sum_{m=1}^k(p_{mj}+p_{jm})\leq \frac{p_j}{5}, \quad j=1,2,\dots,k.
\label{a1a2}
\end{equation}
Note that, in general, it is allowed for all probabilities to depend on the number of summands $n$, that is,
 $X_1,X_2,\dots$ can form triangular arrays (the scheme of series).  Poisson limit occurs when all $p_j=O(n^{-1})$. In our paper, all probabilities are small, though we nevertheless have included the case $p_j=O(1)$. The second assumption essentially reflects requirement for covariance between $X_1$ and $X_2$ to be small.  In \cite{CeVe15},  similar assumptions are made for one-dimensional case.

%%%%%%%%%%%%%%%%%%%%%%%%%%%%%%%%%%%%%%%%%%%%%%%%%%%%%%%%%%%%%%%%%%%%%%%%%%%%%%%%%%%%%%%%%%%%%%%%%%%%%%%%%%%%%%%%%%

Let $\gamma_j=\max(1,\sqrt{np_j})$, $j=1,2,\dots,k$,
\begin{eqnarray*}
\veps_1&:=&\sum_{r,j=1}^k\frac{p_{rj}+p_rp_j}{\gamma_r\gamma_j},\quad
\veps_2:=\sum_{r,j,m=1}^k\frac{p_{rjm}+p_{rj}p_m+p_rp_jp_m}{\gamma_r\gamma_j\gamma_m},\\
\veps_3&:=&
\Ab{\sum_{r,m=1}^k\frac{2p_{rm}-3p_rp_m}{\sqrt{p_rp_m}}},\quad \veps_4:=
\sum_{r,m=1}^k\frac{p_{rm}+p_rp_m}{\sqrt{p_rp_m}}
,\\
\veps_5&:=&\bigg\vert \sum_{r,j,m=1}^k\frac{3p_{rjm}-12p_{rj}p_m+10p_rp_jp_m}{\sqrt{p_rp_jp_m}}
+2\sum_{r=1}^k\frac{3p_{rrr}-12p_{rr}p_r+10p_r^3}{p_r\sqrt{p_r}}\bigg\vert,\\
\veps_6&=&\sum_{r,j,m=1}^k\frac{p_{rjm}+p_rp_{jm}+p_rp_jp_m}{\sqrt{p_rp_mp_j}}.
\end{eqnarray*}

We begin with a result for the Poisson approximation.
%%%%%%%%%%%%%%%%%%%%%%%%%%%%%%%%%%%%%%%%%%%%%%%%%%%%%%%%%%%%%%%%%%%%%%%%%%%%%%
\begin{theorem} \label{P} If the conditions in (\ref{a1a2}) are satisfied, then
 \begin{eqnarray*}
 \norm{F_n-\Pois(\bm{\lambda})}_\infty&\leq& C(k)n\veps_1\prod_{j=1}^k\gamma_j^{-1},\\
 \norm{F_n-\Pois(\bm{\lambda})}_\alpha&\leq& C(k,\alpha)n\veps_1\prod_{j=1}^k\gamma_j^{-(\alpha-1)/\alpha}, \quad  (\alpha\geq 2).
 \end{eqnarray*}
 \end{theorem}
%%%%%%%%%%%%%%%%%%%%%%%%%%%%%%%%%%%%%%%%%%%%%%%%%%%%%%%%%%%%%%%%%%%%%%%%%%%%%%%%%%
Observe that, for $n\veps_1=o(1)$ it suffices that $p_r=o(1), p_{rj}=O(p_rp_j)$, $(r,j=1,\dots,k)$.
Unlike in \cite{roo03,roo17},  Kerstan's method or other convolution technique can not be applied for proofs, since we are dealing with dependent rvs.
The Heinrich method, used in this paper, involves iterations of estimates and results in very large constants. Therefore, we concentrated our efforts on obtaining correct order of estimates leaving the question about the magnitude of constants and their dependence on dimension $k$ for the future research. Note that asymptotic constants can be small, see Proposition \ref{prop} below.
 %in one-dimensional case turned out to be quite small, see \cite{CeVe15}.
As  seen from the following corollary, the order of approximation in Theorem \ref{P} is comparable to known results.

\begin{corollary} If conditions in (\ref{a1a2}) are satisfied and  $np_j\geq 1$, $i=1,2,\dots,k$, then
\begin{equation}
\norm{F_n-\Pois(\bm{\lambda})}_\infty\leq C(k)\prod_{j=1}^k\frac{1}{\sqrt{np_j}}
\bigg(\sum_{r,m=1}^k\frac{p_{rm}}{\sqrt{p_rp_m}}+\bigg(\sum_{r=1}^k\sqrt{p_r}\bigg)^2\bigg).
\label{rez1a}
\end{equation}
\end{corollary}
Observe that for the case $k=1$,   the upper  bound given in (\ref{rez1a})
 coincides with (\ref{ve1}), up to a constant.  Moreover, if $X_1,X_2,\dots$ are independent rvs, then (\ref{rez1a}) is of the same order of accuracy as (\ref{ros2}).
Next we consider  probable improvements  by short asymptotic expansion.

%%%%%%%%%%%%%%%%%%%%%%%%%%%%%%%%%%%%%%%%%%%%%%%%%%%%%%%%%%%%%%%%%%%%%%%%%%%%%%
\begin{theorem}\label{P1} Under the assumptions in  (\ref{a1a2}),  we have
\begin{eqnarray*}
\norm{F_n- \Pois(\bm{\lambda})-A_1}_\infty&\leq& C(k)(n\veps_2+n^2\veps_1^2)\prod_{j=1}^k\gamma_j^{-1},\\
\norm{F_n- \Pois(\bm{\lambda})-A_1}_\alpha&\leq& C(k,\alpha)(n\veps_2+n^2\veps_1^2)\prod_{j=1}^k\gamma_j^{-(\alpha-1)/\alpha},
\quad (\alpha\geq 2).
\end{eqnarray*}
\end{theorem}

It is easy to check that Theorem \ref{P1} is an improvement  over Theorem \ref{P}, if all probabilities $p_r=o(1)$, for $ 1 \leq r \leq k. $ On the other hand,
in the sense of order, there is no difference between corresponding estimates, if all $p_j=O(1)$. Indeed, in the latter case $\veps_1=O(1/n)$ and the accuracy in both cases is of the order $O(n^{-k(\alpha-1)/\alpha})$.
Next consider  approximation with short asymptotic expansion in the exponent, that is signed compound Poisson measure  $G$  defined in Section \ref{intro}.

%%%%%%%%%%%%%%%%%%%%%%%%%%%%%%%%%%%%%%%%%%%%%%%%%%%%%%%%%%%%%%%%%%%%%%%%%%%%%%
\begin{theorem}\label{G}
Let the conditions in  (\ref{a1a2}) hold. Then		
 \begin{eqnarray*}\norm{F_n-G}_\infty&\leq& C(k)n\veps_2\prod_{j=1}^k\gamma_j^{-1},\\
 \norm{F_n-G}_\alpha&\leq& C(k,\alpha)n\veps_2\prod_{j=1}^k\gamma_j^{-(\alpha-1)/\alpha}, \quad (\alpha\geq 2).
 \end{eqnarray*}
 \end{theorem}

%%%%%%%%%%%%%%%%%%%%%%%%%%%%%%%%%%%%%%%%%%%%%%%%%%%%%%%%%%%%%%%%%%%%5
In comparison to Poisson approximation, the signed compound Poisson approximation $G$ is always smaller by the factor $n^{-1/2}$.
Observe also that, when $k=1$, the local  estimate in Theorem \ref{G}  coincides with (\ref{ve2}) up to constant.

\begin{corollary} If conditions (\ref{a1a2}) are satisfied and  $np_j\geq 1$, $i=1,2,\dots,k$, then
\begin{equation*}
\norm{F_n- G}_\infty\leq \frac{C(k)}{\sqrt{n}}\prod_{j=1}^k\frac{1}{\sqrt{np_j}}\sum_{r,j,m=1}^k\frac{p_{rjm}+p_{rj}p_m+p_rp_jp_m}{\sqrt{p_rp_jp_m}}
 .
\label{rez3a}
\end{equation*}
\end{corollary}
%%%%%%%%%%%%%%%%%%%%%%%%%%%%%%%%%%%%%%%%%%%%%%%%%%%%%%%%%%%%%%%%%%%%%%%%%%%%%%%%%%%%%%%%%%%%%%%%%%%%%%%

Are the estimates in Theorems \ref{P} and \ref{G} of the right order? To some extent,
 an affirmative  answer is given by the lower bounds given below.

%%%%%%%%%%%%%%%%%%%%%%%%%%%%%%%%%%%%%%%%%%%%%%%%%%%%%%%%%%%%%%%%%%%%%%%%%%%%%%%%%%%%%%%%%%%%%%%
\begin{theorem} \label{LB} Let the conditions in (\ref{a1a2}) be satisfied and  $np_j\geq 1$ for  $1 \leq j \leq k$, $\alpha\geq 1$.
There exists constants $C_i(k)$, $1 \leq i \leq 6$, such that, for any $b\geq 1$, the following lower bound estimates hold:
\begin{eqnarray*}
\norm{F_n-\Pois(\bm{\lambda})}_\infty&\geq&\frac{C_1(k)}{b^2}\bigg(\veps_3-\frac{C_2(k)}{\min(b,n)}\,\veps_4\bigg)\prod_{j=1}^k(np_j)^{-1/2},\\
\norm{F_n-\Pois(\bm{\lambda})}_\alpha&\geq&\frac{C_3(k)5^{-k(\alpha-1)/\alpha}}{b^2}\bigg(\veps_3-\frac{C_2(k)}{\min(b,n)}\,\veps_4\bigg)
\prod_{j=1}^k(np_j)^{-(\alpha-1)/(2\alpha)},\\
\norm{F_n-G}_\infty&\geq&\frac{C_4(k)}{b^3\sqrt{n}}\bigg(\veps_5-\frac{C_5(k)}{\min(b,n)}\,\veps_6\bigg)\prod_{j=1}^k(np_j)^{-1/2},\\
\norm{F_n-G}_\alpha&\geq&\frac{C_6(k)5^{-k(\alpha-1)/\alpha}}{b^3\sqrt{n}}\bigg(\veps_6-\frac{C_5(k)}{\min(b,n)}\,\veps_6\bigg)
\prod_{j=1}^k(np_j)^{-(\alpha-1)/(2\alpha)}.
\end{eqnarray*}
\end{theorem}

%%%%%%%%%%%%%%%%%%%%%%%%%%%%%%%%%%%%%%%%%%%%%%%%%%%%%%%%%%%%%%%%%%%%%%%%%%%%%%%%%%%%%%%%%%%%%%%%%%%
Note that, unlike Theorems \ref{P}--\ref{G}, we have $\alpha\geq 1$. Therefore, by taking $\alpha=1$, we can establish the lower  estimates for total variation norm.
%%%%%%%%%%%%%%%%%%%%%%%%%%%%%%%%%%%%%%%%%%%%%%%%%%%%%%%%%%%%%%%%%%%%%%%%%%%%%%%%%%%%%%%%%%%%%%%%%%%
%%%%%%%%%%%%%%%%%%%%%%%%%%%%%%%%%%%%%%%%%%%%%%%%%%%%%%%%%%%%%%%%%%%%%%%%%%%%%%%%%%%%%%%%%%%%%%%%%%%
In some cases, estimates in Theorem \ref{LB}  can be trivially negative. Therefore, Theorem \ref{LB}  must be applied   when $\veps_4=O(\veps_3)$. In this case, we can ensure non-triviality of estimates by choosing large enough  constant $b$, so that $C_2(k)/\min(b,n)$ becomes small for all  $n\geq b$. We illustrate this approach by considering independent random vectors. Then $p_{rm}=p_rp_m$ and
\[\veps_3=\bigg(\sum_{r=1}^k\sqrt{p_r}\bigg)^2,\quad \veps_4=2\bigg(\sum_{r=1}^k\sqrt{p_r}\bigg)^2.\]
Therefore, by choosing $b=4C_2(k)$, for all $n\geq 4C_2(k)$ we get
\[\veps_3-\frac{C_2(k)}{\min(b,n)}\,\veps_4\geq \frac{1}{2}\bigg(\sum_{r=1}^k\sqrt{p_r}\bigg)^2.\]
Similarly we can estimate $\veps_5$ and $\veps_6$. We formulate lower estimates for independent vectors as a corollary.
Let as in previous section  $W_n= \tilde X_1+\tilde X_2+\dots+\tilde X_n$ be a sum of independent copies of $X_1$.

%%%%%%%%%%%%%%%%%%%%%%%%%%%%%%%%%%%%%%%%%%%%%%%%%%%%%%%%%%%%%%%%%%%%%%%%%%%%%%%%%%%%%%%%%%%%%%%%%%%

\begin{corollary} \label{LBIND} Let $\max_{1\leq j\leq k}p_j\leq 1/(144k)$  and  $np_j\geq 1$, $i=1,2,\dots,k$, $\alpha\geq 1$.
Then there exist absolute constants $C_{7}(k), C_{8}(k), C_9(k)$,  such that, for any $n\geq C_{7}(k)$, the following estimates hold
\begin{eqnarray*}
\norm{\eL(W_n)-\Pois(\bm{\lambda})}&\geq& C_{8}(k)\bigg(\sum_{r=1}^k\sqrt{p_r}\bigg)^2 \prod_{j=1}^k(np_j)^{-1/2},\\
\norm{\eL(W_n)-\Pois(\bm{\lambda})}_\alpha&\geq& C_{8}(k)5^{-k(\alpha-1)/\alpha}\bigg(\sum_{r=1}^k\sqrt{p_r}\bigg)^2
\prod_{j=1}^k(np_j)^{-(\alpha-1)/(2\alpha)}, \\
\norm{\eL(W_n)-G}_\infty&\geq&\frac{C_{9}(k)p^{3/2}}{\sqrt{n}}\prod_{j=1}^k(np_j)^{-1/2},\\
\norm{\eL(W_n)-G}_\alpha&\geq&\frac{C_{9}(k)p^{3/2}}{\sqrt{n}}
\prod_{j=1}^k(np_j)^{-(\alpha-1)/(2\alpha)}.
\end{eqnarray*}
\end{corollary}

%%%%%%%%%%%%%%%%%%%%%%%%%%%%%%%%%%%%%%%%%%%%%%%%%%%%%%%%%%%%%%%%%%%%%%%%%%%%%%%%%%%%%%%%%%%%%%%%%%%
Comparing Corollary \ref{LBIND} with (\ref{ros2}) and (\ref{ZHwang}) we see that lower estimates have the same order as the  upper estimates. The same reasoning applies for dependent vectors.
For example, if $p_rp_m=o(p_{rm})$, then we always can choose $b$ in such a way that, for sufficiently large $n$,
\[\veps_3-\frac{C_2(k)}{\min(b,n)}\,\veps_4\geq C(k) \sum_{r,m=1}^k\frac{p_{rm}}{\sqrt{p_rp_m}}.
\]
Similarly, we can ensure  the first right-hand side estimate in Theorem \ref{LB} to be positive if $p_{mr}=o(p_mp_r)$. Application of Theorem \ref{LB} to one well-known statistic is given below.
%%%%%%%%%%%%%%%%%%%%%%%%%%%%%%%%%%%%%%%%%%%%%%%%%%%%%%%%%%%%%%%%%%%%%%%%%%%%%%%%%%%%%%%%%%%%%%%%%%%%%%%%%%%%%%%%%%%%%%%%

%%%%%%%%%%%%%%%%%%%%%%%%%%%%%%%%%%%%%%%%%%%%%%%%%%%%%%%%%%%%%%%%%%%%%%%%%%%%%%%%%%%%%%%%%%%%%%%%%%%%%%%%%%%%%%%%%%%%%%%%

%%%%%%%%%%%%%%%%%%%%%%%%%%%%%%%%%%%%%%%%%%%%%%%%%%%%%%%%%%%%%%%%%%%%%%%%%%%%%%%%%%%%%%%%%%%%%%%%%%%%%%%%%%%%%%%%%%%%%%%%

\section{An Application} Let $\xi_1,\xi_2,\dots$ be independent identically distributed Bernoulli variables with probability $q\in (0,1)$.
One of the best known and thoroughly investigated examples of the sum of 1-dependent random variables is 2-runs statistic, that is, $S=\xi_1\xi_2+\xi_2\xi_3+\cdots+\xi_n\xi_{n+1}$, see \cite{CeVe15,PeCe10,WXia08} and the references therein. We similarly construct 2-dimensional  parallel runs with  random switching between them. More precisely, let $\bar{\xi}_1, \bar{\xi}_2,\dots$ be a sequence of independent
 identically distributed Bernoulli random variables with probability $\bar{q}\in(0,1)$  and let $\eta_1,\eta_2,\dots$ be another sequence of independent identically distributed  Bernoulli variables with probability $\delta\in(0,1)$. Moreover, we assume the random variables in all three sequences to be mutually independent. Let us define a sequence of 1-dependent 2-dimensional rvs in the following way:
 \[Y_j=(\eta_j\xi_j\xi_{j+1},(1-\eta_j)\bar{\xi}_j\bar{\xi}_{j+1}),\quad j=1,2,\dots.\]

 Let  $S=Y_1+Y_2+\cdots+Y_n$. It is easy to check that $p_1=\delta q^2$, $p_2=(1-\delta)\bar{q}^2$, $p_{11}=\delta^2q^3$, $p_{12}=p_{21}=\delta(1-\delta)q^2\bar{q}^2$, $p_{22}=(1-\delta)^2\bar{q}^3$, $p_{111}=\delta^3q^4$,
 $p_{112}=p_{211}=\delta^2(1-\delta)q^3\bar{q}^2$, $p_{121}=\delta^2(1-\delta)q^4\bar{q}^2$, $p_{122}=p_{221}=\delta (1-\delta)^2q^2\bar{q}^3$,
 $p_{212}=\delta(1-\delta)^2q^2\bar{q}^4$, $p_{222}=(1-\delta)^3\bar{q}^4$.

 Let us assume that $q,\bar{q}\leq 1/17$. Then  conditions (\ref{a1a2}) are satisfied.
Observe also that
 \[\sum_{r,m=1}^2\frac{p_{rm}}{\sqrt{p_rp_m}}=\delta q+2\sqrt{\delta(1-\delta)}q\bar{q}+(1-\delta)\bar{q}\leq 2(q\delta +\bar{q}(1-\delta)).\]
 Therefore, when $np_1,np_2\geq 1$,  it follows from (\ref{rez1a}) and Theorem \ref{LB} that
 there exists a constant  $C_{10}$ such that for $n\geq C_{10}$,
 \begin{equation}
\frac{C_{11}(q\delta +\bar{q}(1-\delta))}{nq\bar{q}\sqrt{\delta(1-\delta)}}\leq  \norm{\eL(S)-\Pois(\bm{\lambda})}_\infty\leq
\frac{C_{12}(q\delta +\bar{q}(1-\delta))}{nq\bar{q}\sqrt{\delta(1-\delta)}},
\label{exP}
\end{equation}
where $\bm{\lambda}=(n\delta q^2,n(1-\delta)\bar{q}^2)$.
Similarly, for the case $np_1,np_2\geq 1$ there exists $C_{13}$ such that for all $n\geq C_{13}$, we obtain
\begin{eqnarray}
\frac{C_{14}(q\delta\sqrt{\delta}+\bar{q}(1-\delta)\sqrt{1-\delta}\,)}{n\sqrt{n}\,q\bar{q}\sqrt{\delta(1-\delta)}}&\leq& \norm{\eL(S)- G}_\infty
\nonumber\\
&\leq&
\frac{C_{15}(q\delta\sqrt{\delta}+\bar{q}(1-\delta)\sqrt{1-\delta}\,)}{n\sqrt{n}\,q\bar{q}\sqrt{\delta(1-\delta)}}.
\label{exG}
\end{eqnarray}
The condition for $n$ to be larger than some absolute constant (which can be estimated with the help of Lemma \ref{lblem} below)  is needed for lower bounds only. The upper bounds  in (\ref{exP}) and (\ref{exG}) hold for all $n\geq 1$. As expected, the benefits of expansion in the exponent are expressed through additional factor $1/\sqrt{n}$.
In one dimensional case, the local estimate for $G$ is $Cn^{-1}$, see \cite{PeCe10}, Theorem 2. The additional multiplier $1/\sqrt{nq\bar{q}}$ appears because we investigate two-dimensional vectors.
When $q$ and $\bar{q}$ are slowly vanishing, the explicit form of the rvs allows to estimate asymptotic constant.

\begin{proposition} \label{prop} Let $\delta$ be a constant, $\max(q,\bar{q})=o(1)$, and $\min(nq,n\bar{q})\to \infty$, as $n\to\infty$. Then
\begin{align*}
& \lim_{n\to\infty} \frac{nq\bar{q}\sqrt{\delta(1-\delta)}}{(q\delta +\bar{q}(1-\delta))}
\norm{\eL(S)-\Pois(\bm{\lambda})}_\infty\\
 & \hspace{2cm} \leq\frac{1}{\ee}\bigg(1+\sqrt{\frac{\pi}{2}}\bigg)\bigg\{\frac{1}{\sqrt{6}}\bigg(1+\sqrt{\frac{\pi}{4}}\bigg)+
\frac{1}{8}\bigg(1+\sqrt{\frac{\pi}{2}}\bigg)\bigg\}\\
&  \hspace{2cm} =0.871\dots~ .\end{align*}
\end{proposition}

Proposition \ref{prop} serves as an indicator that constants in above theorems should not be very large.
%Similarly, we can check other cases and norms. For example, if
%$np_1\geq 1$, $np_2\leq 1$, then
%\[\norm{\eL(S)-\Pois(\bm{\lambda})}_\infty\leq C\bigg(\frac{\delta p}{n}+(1-\delta)^2\bar{p}^3\bigg),\]
%etc.
%%%%%%%%%%%%%%%%%%%%%%%%%%%%%%%%%%%%%%%%%%%%%%%%%%%%%%%%%%%%%%%%%%%%%%%%%%%%%%%%%%%%%%%%%%%%%%%%%%%%%%%%%%%%%%%%%%%
%AUXILIAR

\section{Auxiliary results}

We  begin from  relating Fourier transforms to local and $\ell_2$  norms.
%%%%%%%%%%%%%%%%%%%%%%%%%%%%%%%%%%%%%%%%%
%%%%%%%%%%%%%%%%%%%%%%%%%%%%%%%%%%%%%%%%%%%%%%%%%%%%%%%%%%%%%%%%%%%%%%%%%%%%%%%%%%%%%%%%%%%%%%%%%%%%%%%%%%%%%%%%%%%%%%%%
%FInversion

\begin{lemma} \label{locinv} Let $M$ be (signed) measure concentrated on $\ZZ^k$.
Then
\begin{eqnarray*}
\norm{M}_\infty&\leq&\frac{1}{(2\pi)^k}\int_{-\pi}^\pi\cdots\int_{-\pi}^\pi\ab{\w
M(\boldt)}\dd t_1\dots\dd t_k, \label{apvertimas}\\
\norm{M}_2&=&\frac{1}{(2\pi)^k}\int_{-\pi}^\pi\cdots\int_{-\pi}^\pi \ab{\w M(\boldt)}^2\dd t_1\dots\dd t_k.
\label{el2inv}
\end{eqnarray*}
\end{lemma}

\emph{Proof.}
 The first inequality follows directly from the inversion formula. The second is multidimensional Parseval's identity.
 However, in order to keep the paper self-contained, we give an outline of the proof. First we introduce measure $M^{-}\{\boldm\}=M\{-\boldm\}$. One one hand, convolution of both measures at point zero is equal to
\[M*M^{-}\{\boldzero\}=\sum_{\boldm\in\ZZ^k}M^2\{\boldm\}.\]
On the other hand, by inversion formula
\[M*M^{-}\{\boldzero\}=\frac{1}{(2\pi)^k}\int_{-\pi}^\pi\cdots\int_{-\pi}^\pi \w M(\boldt)\w M(-\boldt)\dd t_1\dots\dd t_k.\qquad
\hbox{\qed}
\]

%%%%%%%%%%%%%%%%%%%%%%%%%%%%%%%%%%%%%%%%%%%%%%%%%%%%%%%%%%%%%%%%%%%%%%%%%%%%%%%%%%%%%%%%%%
For the lower bounds, an  appropriate inversion formula is needed. First, we introduce additional notation. Let
\begin{equation*}
\bolda:=(a_1,a_2,\dots,a_k),\quad \boldbeta:=(\beta_1,\beta_2,\dots,\beta_k),
\quad \boldt_{\boldbeta} :=\Big(\frac{t_1}{\beta_1},\frac{t_2}{\beta_2},\dots,\frac{t_k}{\beta_k}\Big).
\label{lb1}
\end{equation*}
Next we define some weight functions. Lemma below holds true if we
\[\psi_j(t_j)=\ee^{-t_j^2/2}\qquad\hbox{or}\qquad \psi_j(t_j)=t_j\ee^{-t^2_j/2}.\]

 \begin{lemma} \label{lblem}  Let  $M$ be a finite measure concentrated on $\ZZ^k$. Then, for any
  $\bolda\in\RR^k$ and $\beta_j\geq 1$, $(j=1,2,\dots,k)$, the following inequalities hold:
  \begin{eqnarray*}
  \norm{M}_\infty&\geq& (4\sqrt{2\pi})^{-k}\prod_{j=1}^k\beta_j^{-1}
  \ab{V(\bolda,\boldbeta)},\quad \norm{M}\geq (\sqrt{2\pi})^{-k} \ab{V(\bolda,\boldbeta)},
  \label{lb2}\\
  \norm{M}_\alpha&\geq&(\sqrt{2\pi})^{-k}5^{-k(\alpha-1)/\alpha}\prod_{j=1}^k \beta_j^{-(\alpha-1)/\alpha}
  \ab{V(\bolda,\boldbeta)},\quad (\alpha> 1).
  \end{eqnarray*}
  Here
  \[
    V(\bolda,\boldbeta):=\int_{-\infty}^\infty\cdots\int_{-\infty}^\infty  \prod_{m=1}^k\psi_m(t_m)\ee^{-\ii(\boldt_{\boldbeta},\bolda)}
  \w M(\boldt_{\boldbeta})\dd t_1\dd t_2\dots\dd t_k.\]
  \end{lemma}

 \emph{Proof.} We adopt the proof of Lemma 10.1 from \cite{Ce16}. Observe that
 \[\ee^{-\ii(\boldt_{\boldbeta},\bolda)}  \w M(\boldt_{\boldbeta})=\sum_{\boldm\in\ZZ^k}\ee^{\ii(\boldt_{\boldbeta},\boldm-\bolda)}M\{\boldm\}.\]
   By interchanging the order of
 integration and summation, we obtain
 \begin{eqnarray*}
 V(\bolda,\boldbeta)&=&\sum_{\boldm\in\ZZ^k}\prod_{j=1}^k\bigg(\int_{-\infty}^\infty \psi_j(t_j)\ee^{\ii t_j(m_j-a_j)/\beta_j}\dd t_j\bigg)M\{\boldm\}\nonumber\\
 &=&(\sqrt{2\pi})^k \sum_{\boldm\in\ZZ^k}\prod_{j=1}^k\tilde\psi_j\Big(\frac{m_j-a_j}{\beta_j}\Big)M\{\boldm\}.
 \label{lbp}
 \end{eqnarray*}
 Here, depending on the choice of $\psi_j(t_j)$,
 \[\tilde\psi_j(y)=\ee^{-y^2/2} \quad\hbox{or}\quad \tilde\psi_j(y)=\ii y\ee^{-y^2/2}.\]
  For the norm $\ell_\alpha$, $\alpha >1$,  we apply H\"older's inequality
 \begin{eqnarray*}
 \ab{V(\bolda,\boldbeta)}&\leq& (\sqrt{2\pi})^k\norm{M}_\alpha\bigg(
 \sum_{\boldm\in\ZZ^k}\prod_{j=1}^k\tilde\psi_j^{\alpha/(\alpha-1)}\Big(\frac{m_j-a_j}{\beta_j}\Big)\bigg)^{(\alpha-1)/\alpha}\\
 &=& (\sqrt{2\pi})^k\norm{M}_\alpha \bigg(\sum_{m_j\in\ZZ}\tilde\psi_j\Big(\frac{m_j-a_j}{\beta_j}\Big)\bigg)^{(\alpha-1)/\alpha}.
 \end{eqnarray*}
 Let $q:=\alpha/(\alpha-1)$, $y_j:=(m_j-a)/\beta_j$. Then, since $q>1$,
 \[\sum_{j\in\ZZ}\ee^{-qy^2_j/2}\leq 1+\beta_j\frac{\sqrt{2\pi}}{\sqrt{q}}\leq 1+\beta_j\sqrt{2\pi}\]
 and
 \begin{eqnarray*}
 \sum_{j\in\ZZ}\ab{y_j}^q\ee^{-qy^2_j/2}&\leq&\frac{1}{2}\sum_{j\in\ZZ}(1+y_j^{2q})\ee^{-qy_j^2/\ee}\ee^{-q(\ee-2)y_j^2/(2\ee)}\\
 &\leq&\frac{1}{2}\Big(1+\beta_j\sqrt{2\pi}+\max_{x>0}x^q\ee^{-qx/\ee}\sum_{j\in\ZZ}\ee^{-q(\ee-2)y_j^2/(2\ee)}\Big)\\
 &\leq&\frac{1}{2}\Big(1+\beta_j\sqrt{2\pi}+1+\beta_j\sqrt{2\pi\ee/(\ee-2)}\Big)<5\beta_j.
 \end{eqnarray*}
 Similarly, for the local norm,
 \begin{eqnarray*}
 \ab{V(\bolda,\boldbeta)}&\leq& (\sqrt{2\pi})^k\norm{M}_\infty\sum_{\boldm\in\ZZ^k}\prod_{j=1}^k\tilde\psi_j\Big(\frac{m_j-a_j}{\beta_j}\Big)\\
 &=& (\sqrt{2\pi})^k\norm{M}_\infty\prod_{j=1}^k\bigg(\sum_{m_j\in\ZZ}\tilde\psi_j\Big(\frac{m_j-a_j}{\beta_j}\Big)
 \bigg)\\
 & \leq& (\sqrt{2\pi})^k\norm{M}_\infty\prod_{j=1}^k(4\beta_j).%\qquad\qquad\qed
 \end{eqnarray*}
  For the total variation norm
 \[
 \ab{V(\bolda,\boldbeta)}\leq (\sqrt{2\pi})^k \sum_{\boldm\in\ZZ^k}\prod_{j=1}^k
 \Ab{\tilde\psi_j\Big(\frac{m_j-a_j}{\beta_j}\Big)}\ab{M\{\boldm\}}\leq
  (\sqrt{2\pi})^k\norm{M}.\qquad\hbox{\qed}\]
 %   \qed

 %%%%%%%%%%%%%%%%%%%%%%%%%%%%%%%%%%%%%%%%%%%%%%%%%%%%%%%%%%%%%%%%%%%%%%%%%%%%%%%%%%%%%%%%%%%%%%%

 Next we formulate two technical results.

%%%%%%%%%

%%%%%%%%%%%%%%%%%%%%%%%%%%%%%%%%%%%%%%%%%%%%%
\begin{lemma}\label{shorgin} Let $a>0$, $m\geq 1$. Then
\begin{eqnarray}\frac{1}{2\pi}\int_{-\pi}^\pi\ab{\sin(t/2)}^m\exponent{-a\sin^2(t/2)}\dd
t&\leq&\sqrt{\ee}\bigg(1+\sqrt{\frac{\pi}{2}}\,\bigg)\bigg(\frac{m}{2a\ee}\bigg)^{(m+1)/2},
\label{sh1}\\
\frac{1}{2\pi}\int_{-\pi}^\pi\exponent{-a\sin^2(t/2)}\dd
t&\leq&\bigg(1+\sqrt{\frac{\pi}{2}}\,\bigg)\frac{1}{\sqrt{6a}}.
\label{sh2}
\end{eqnarray}
\end{lemma}
Lemma \ref{shorgin} is essentially
Lemma 6 in \cite{roo01}.

%%%%%%%%%%%%%%%%%%%%%%%%%%%%%%%%%%%%%%%%%%%%%%%%%%%%%%%%%%%%%%%%%%%%%%%%%%%%%%%%

%ReF lemma

\begin{lemma}  Let $0\leq p\leq 1$, $a^2+b^2\leq 1$.
Then
\begin{equation}
\ab{(1-p)+p(a+\ii b)}\leq 1+p(1-p)(a-1)\leq\exponent{p(1-p)(a-1)}. \label{chf0}
\end{equation}
\end{lemma}
The proof is trivial and can be found, for example, in \cite{KC14}.

%%%%%%%%%%%%%%%%%%%%%%%%%%%%%%%%%%%%%%%%%%%%%%%%%%%%%%%%%%%%%%%%%%%%%%%%%%%%%%%%%%%%%%%%%%%%%%%%%%%%%%%%%%%%%%%%%%%%%%%%%
%%%%%%%%%%%%%%%%%%%%%%%%%%%%%%%%%%%%%%%%%%%%%%%%%%%%%%%%%%%%%%%%%%%%%%%%%%%%%%%%%%%%%%%%%%%%%%%%%%%%%%%%%%%%%%%%%
For $X_j=(X_{j1},X_{j2},\dots,X_{jk})$, we introduce accompanying complex-valued random variables
\begin{equation*}
Z_j:=\exponent{\ii(\boldt, X_j)}-1,\quad (\boldt,X_j)=t_1X_{j1}+t_2X_{j2}+\cdots+t_kX_{jk}.
\label{aux1}
\end{equation*}
We assume that $\w\Expect (Z_1)=\Expect Z_1$, $\wE(Z_1,Z_2)=\EE Z_1Z_2-\EE Z_1\EE Z_2$ and, for
$m\geq 3$, define
\begin{equation*}
 \w\Expect (Z_1,Z_2,\cdots, Z_m)=\Expect Z_1Z_2\cdots
Z_k-\sum_{j=1}^{m-1}\w\Expect (Z_1,\cdots ,Z_j)\Expect
Z_{j+1}\cdots Z_{m}. \label{capY}
\end{equation*}

%%%%%%%%%%%%%%%%%%%%%%%%%%%%%%%%%%%%%%%%%%%%%%%%%%%%%%%%%%
%%%%%%%
%%%%%%%%%%%%%%%%%%%%%%%%%%%%%%%%%%%%%%%%%%%%%%%%%%%%%%%%%%%%%%%%%%%%%%%%%%%%%%%%%%%%%%%%%%%%%%%%%%%%%%%%%%%%%%%%%%%%%%%%%
The essence of Heinrich's method is the following characterization lemma.

\begin{lemma}\label{Hversion} Let $Z_1,\dots,Z_n$ be defined as  above and let
\begin{equation}
\sqrt{\Expect\ab{Z_j}^2}\leq \frac{1}{6},\qquad j=1,2,\dots,n.\label{Hcond}
\end{equation}
 Then
\[\w F_n(\boldt)=\vfi_1(\boldt)\vfi_2(\boldt)\dots\vfi_n(\boldt),\]
 where $\vfi_1(\boldt)=\Expect Z_1$ and, for $m=2,\dots, n$,
\[\vfi_m(\boldt)=1+ \Expect Z_m+\sum_{j=1}^{m-1}\frac{\wE(Z_j,Z_{j+1},\dots,Z_m
)}{\vfi_j(\boldt)\vfi_{j+1}(\boldt)\dots \vfi_{m-1}(\boldt)}.\]
\end{lemma}

Lemma \ref{Hversion} follows from more general Lemma
3.1  in \cite{H82}, see also Theorem 1 in \cite{H85}.
Also, the next lemma also can be found in \cite{H82}.
%%%%%%%%%%%%%%%%%%%%%%%%%%%%%%%%%%%%%%%%%%%%%%%%%%%%%%%%%%%%%%%%%%%%%%%%%%%%%%%%%%%%%%%%%%%%%%%%%%%%%%%%%%%%%%%%%%%%%%%%%
%%%\WE HOLDER
%%%%%%%%%%%%%%%%%%%%%%%%%%%%%%%%%%%%%%%%%%%%%%%%%%%%%%%%%%%%%%%%%%%%%%%%%%%%%%%%%%%%%%%%%%%%%%%%%%%%%%%%%%%%%%%%%%
\begin{lemma} \label{Hei3aa} Let $Z_1,Z_2,\dots,Z_k$ be 1-dependent complex-valued random variables with
$\Expect\ab{Z_j}^2<\infty$,  $1 \leq j \leq m. $ Then
\[
\ab{\wE (Z_1, Z_2, \cdots, Z_k)}\leq
2^{k-1}\prod_{j=1}^k(\Expect\ab{Z_j}^2)^{1/2}.
\]
\end{lemma}
%%%%%%%%%%%%
For the sake of convenience, we collect all  the facts about $Z_j$ and present it  in the following lemma.  Let
\[u(\boldt):=\sum_{r=1}^kp_r \sin^2(t_r/2).\]

\begin{lemma}\label{Zj} Let the assumptions in (\ref{a1a2}) hold. Then
\begin{eqnarray}
 \EE Z_1&=&\sum_{r=1}^kp_r(\ee^{\ii t_r}-1),\quad
\EE Z_1Z_2=\sum_{r,m=1}^kp_{rm}(\ee^{\ii t_r}-1)(\ee^{\ii t_m}-1),\label{z0}\\
  \EE \ab{Z_1}&=&2\sum_{r=1}^kp_r\ab{\sin(t_r/2)},\quad \EE\ab{Z_1}^2=4u(\boldt),\quad \ab{Z_j}\leq2,
\label{z1}\\
 \EE\ab{Z_1Z_2}&\leq& 0.4u(\boldt),\quad \EE\ab{Z_1}\EE\ab{Z_2}\leq
\frac{u(\boldt)}{36},\label{z2}\\
 Re\,\EE Z_1&=&-2u(\boldt), \quad \ab{\wE(Z_j,Z_{j+1},\dots,Z_m)}\leq 8 u(\boldt)(4u(\boldt))^{m-j-1}.
\label{z3}
\end{eqnarray}
Here $Re\,\EE Z_1$ denotes the real part of complex number $\EE Z_1$.
\end{lemma}
%%%%%%%%%%%%%%%%%%%%%%%%%%%%%%%%%%%%%%%%%%%%%%%%%%%%%%%%%%%%%%%%%%%%%%%%%%%%%%%%%%%%%%%%%%%%%%%%%%%%%%%%%%%%%%%%%%%%%%%%%

\emph{Proof}. Let $ a+ib$ be a complex number. Clearly,  $\ab{a+\ii b}^2=a^2+b^2$. Therefore, $\ab{Z_j}\leq \ab{\cos(\boldt,X_j)
+\ii\sin(\boldt,X_j)}+1\leq 2$. Observe that due to the first assumption in (\ref{a1a2})
\begin{eqnarray*}
\EE\ab{Z_1}\EE\ab{Z_2}&\leq& 4\sum_{r,m=1}^kp_rp_m\ab{\sin(t_r/2)}\ab{\sin(t_m/2)}\\
&\leq&
2\sum_{r,m=1}^k(p_2^2\sin^2(t_r/2)+p_m^2\sin^2(t_m/2))\\
&=&4k\sum_{j=1}^k p_j^2\sin^2(t_j/2)\leq \frac{u(\boldt)}{36}.
\end{eqnarray*}
Similarly, due to the second assumption in (\ref{a1a2}),
\begin{eqnarray}
\EE\ab{Z_1Z_2}&=&4\sum_{r,m=1}^kp_{rm}\ab{\sin(t_r/2)}\ab{\sin(t_m/2)}\nonumber\\
&\leq& 2\sum_{r,m=1}^kp_{rm}(\sin^2(t_r/2)+\sin^2(t_m/2))\nonumber\\
&=&2\sum_{r=1}^k\sin^2(t_r/2)\sum_{m=1}^k(p_{rm}+p_{mr})\leq 0.4\sum_{r=1}^kp_r\sin^2(t_r/2).
\label{vai}\end{eqnarray}
Finally, by Lemma \ref{Hei3aa}
\[ \ab{\wE(Z_j,Z_{j+1},\dots,Z_m)}\leq 2^{m-j}(\sqrt{\EE\ab{Z_1}^2})^{m-j+1}=2^{m-j}\big(2\sqrt{u(\boldt)}\big)^{m-j+1}.\]
%%%%%%%%%%%%%%%%%%%%%%%%%%%%%%%%%%%%%%%%%%%%%%%%%%%%%%%%%%%%%%%%%%%%%%%%%%%%%%%%%%%%%%%%%%%%%%%%%%%%%%%%%%%%%%%%%%%%%%%%%
All other relations follow directly from the definition of $Z_j$.
\qed

%%%%%%%%%%%%%%%%%%%%%%%%%%%%%%%%%%%%%%%%%%%%%%%%%%%%%%%%%%%%%%%%%%%%%%%%%%%%%%%%%%%%%%%%%%%555
\begin{lemma}\label{C1} Assume the conditions in (\ref{a1a2}) hold
	and $\ab{t_j} \leq \pi$, for  $ 1 \leq j \leq k .$
\begin{equation}
\ab{\vfi_m(\boldt)-1}\leq\frac{1}{10},\quad\frac{1}{\ab{\vfi_m(\boldt)}}\leq\frac{10}{9},\quad
\ab{\vfi_m(\boldt)-1-\EE Z_m}\leq 1.93\, u(\boldt),
\label{C2}
\end{equation}
for $ 1 \leq m \leq n.$
\end{lemma}

\emph{Proof.} Further in the proofs, for the sake of brevity, we write $\vfi_m$ instead of $\vfi_m(\boldt)$ whenever no ambiguity can arise. The first two estimates follow from the third one. Indeed, due to (\ref{a1a2}), (\ref{z1}) and the definition of $u(\boldt)$,
\[\ab{\vfi_m-1}\leq\EE\ab{Z_m}+\ab{\vfi_m-1-\EE Z_m}\leq 4\sum_{r=1}^kp_r\leq\frac{4}{144}<\frac{1}{10}.\]
Similarly,
\[\ab{\vfi_m}=\ab{1-\vfi_m-1}\geq 1-\ab{\vfi_m-1}\geq \frac{9}{10}.\]

The proof of the third estimate in (\ref{C2}) is done by mathematical induction. Observe that, due to (\ref{a1a2}) and (\ref{z1}),
condition (\ref{Hcond}) is satisfied and we can apply Lemma \ref{Hversion}. Let all the estimates hold for $j=1,2,\dots,m-1$, $m>4$. (For $m=2,3,4$ the proof is similar and shorter). Then
\begin{eqnarray}
\ab{\vfi_m-1-\EE Z_m}&\leq& \frac{10}{9}\ab{\wE(Z_1,Z_2)}+\bigg(\frac{10}{9}\bigg)^2\ab{\wE(Z_1,Z_2,Z_3)}\nonumber\\
&&+\bigg(\frac{10}{9}\bigg)^3\ab{\wE(Z_1,Z_2,Z_3,Z_4)}\nonumber\\
&&+\sum_{j=1}^{m-5}\bigg(\frac{10}{9}\bigg)^{k-j}\ab{\wE(Z_j,\dots,Z_m)}.
\label{C3}
\end{eqnarray}
By (\ref{z3})
\begin{eqnarray*}
\lefteqn{\sum_{j=1}^{m-5}\bigg(\frac{10}{9}\bigg)^{k-j}\ab{\wE(Z_j,\dots,Z_m)}
\leq \sum_{j=1}^{m-5}\bigg(\frac{10}{9}\bigg)^{k-j}8u(\boldt)(4\sqrt{u(\boldt)})^{k-j-1}}\hskip 2.5cm\\
&\leq&8u(\boldt)\sum_{j=1}^{m-5}\bigg(\frac{10}{9}\bigg)^{k-j}\bigg(\frac{1}{3}\bigg)^{k-j-1}\leq 0.266 u(\boldt).
\end{eqnarray*}
By Lemma \ref{Zj}
\begin{eqnarray*}
\EE\ab{Z_1}&\leq&\frac{1}{72},\quad \ab{u(\boldt)}\leq\frac{1}{144},\quad \EE\ab{Z_1Z_2Z_3}\leq 2\EE\ab{Z_1Z_2}\leq 0.8u(\boldt),\\
\EE\ab{Z_1Z_2Z_3Z_4}&\leq&\sqrt{\EE\ab{Z_1Z_3}\EE\ab{Z_2Z_4}}=(4u(\boldt))^4\leq\frac{16 u(\boldt)}{144}=\frac{u(\boldt)}{9},\\
\ab{\wE(Z_1,Z_2)}&\leq&\EE\ab{Z_1Z_2}+\EE\ab{Z_1}\EE\ab{Z_2}\leq 0.4u(\boldt)+\frac{u(\boldt)}{36}\leq 0.428u(\boldt),\\
\ab{\wE(Z_1,Z_2,Z_3)}&\leq&\EE\ab{Z_1Z_2Z_3}+\ab{\wE(Z_1,Z_2)}\EE\ab{Z_1}+\EE\ab{Z_1}\EE\ab{Z_2Z_3}\leq  0.812 u(\boldt).
\end{eqnarray*}
Similarly we prove that $\ab{\wE(Z_1,Z_2,Z_3,Z_4)}\leq 0.135u(\boldt)$. Substituting these estimates into (\ref{C3}), we complete lemma's proof. \qed

%%%%%%%%%%%%%%%%%%%%%%%%%%%%%%%%%%%%%%%%%%%%%%%%%%%%%%%%%%%%%%%%%%%%%%%%%%%%%%%%%%%%%%%%%%%%%%%
\begin{lemma} \label{expofi} Let the conditions stated in (\ref{a1a2}) be satisfied. Then,
\begin{equation*}
\ab{\vfi_m(\boldt)}\leq 1-0.05u(\boldt)\leq \exponent{-0,05u(\boldt)},\quad
\ab{\w F_n(\boldt)}\leq\exponent{-0.05nu(\boldt)}, \label{chf}
\end{equation*}
for $ 1 \leq m \leq n.$
\end{lemma}

\emph{Proof.} We have
\[\ab{\vfi_m(\boldt)}\leq \ab{1+\EE Z_m+\vfi_m(\boldt)-1-\EE Z_m}\leq \ab{1+\EE Z_m}+1.93u(\boldt).
\]
It remains to apply (\ref{chf0}) to the first summand.
The second estimate in (\ref{chf}) follows  from the Heinrich's decomposition $\w F_n(\boldt)=\vfi_1(\boldt)\cdots\vfi_n(\boldt)$.
\qed

%To make notation shorter below we assume that $\sum_{j=m}^s=0$, if $m>s$.
%%%%%%%%%%%%%%%%%%%%%%%%%%%%%%%%%%%%%%%%%%%%%%%%%%%%%%%%
Let us now denote the remainder terms by
\begin{eqnarray*}
r_1(\boldt)&:=&\sum_{m,j=1}^k(p_{mj}+p_mp_j)\ab{\sin(t_m/2)}\ab{\sin(t_j/2)},\label{r1}\\
r_2(\boldt)&:=&\sum_{i=1}^k\ab{t_i}\sum_{m,j=1}^k(p_{mj}+p_mp_j)\ab{\sin(t_m/2)}\ab{\sin(t_j)},\label{r2}\\
r_3(\boldt)&:=&\sum_{l,m,j=1}^k(p_{lmj}+p_{lm}p_j+p_lp_mp_j)\ab{\sin(t_l/2)}\ab{\sin(t_m/2)}\ab{\sin(t_j/2)},\label{r3}\\
r_4(\boldt)&:=& \sum_{i=1}^k\ab{t_i}\sum_{l,m,j=1}^k(p_{lmj}+p_{lm}p_j+p_lp_mp_j)\ab{t_l}\ab{t_m}\ab{t_j}.\label{r4}
\end{eqnarray*}

\begin{lemma} \label{vfiskl} Let conditions (\ref{a1a2}) be satisfied. Then, for all $r=2,\dots,n$ % and $\ab{t_1},\dota,\ab{t_k}\leq\pi$
\begin{eqnarray*}
\vfi_r(\boldt)&=&1+\EE Z_r+\theta C(k)r_1(\boldt),\\
\vfi_r(\boldt)&=&1+\EE Z_r-\sum_{j,m=1}^k(p_{jm}-p_jp_m)t_jt_m+\theta C(k)r_2(\boldt),\\
\vfi_r(\boldt)&=&1+\EE Z_r+\wE(Z_r,Z_{r-1})
%\sum_{j=1}^kp_j(\ee^{\ii t_j}-1)+\sum_{j,m=1}^k(p_{jm}-p_jp_m)(\ee^{\ii t_j}-1)(\ee^{\ii t_m}-1)
+\theta C(k)r_3(\boldt),
\end{eqnarray*}
and for $r=3,4,\dots,n$,
\begin{eqnarray*}\vfi_r(\boldt)&=&1+\EE Z_r+\wE(Z_r,Z_{r-1})\\
&&+\ii^3\sum_{l,m,j=1}^k(p_{lmj}-3p_{lj}p_m+2p_lp_jp_m)t_lt_mt_j+\theta C(k)r_4(\boldt).
\end{eqnarray*}
\end{lemma}
\emph{Proof.}
 %As before, we write $\vfi_r$ instead of $\vfi_r(\boldt)$.
Let us assume that $k\geq 5$. Applying Lemma \ref{Hversion} we obtain
\begin{equation*}
\vfi_r=1+\EE Z_r+\wE(Z_r,Z_{r-1})
\sum_{j=r-4}^{r-1}\frac{\wE(Z_j,\dots,Z_r)}{\vfi_j\cdots\vfi_{r-1}}
+\sum_{j=1}^{r-5}\frac{\wE(Z_j,\dots,Z_r)}{\vfi_j\cdots\vfi_{r-1}}.
\label{fi1}
\end{equation*}
Estimating the absolute value of the second expression, as in the proof of Lemma \ref{C1}, we get
\[\sum_{j=1}^{r-5}\frac{\ab{\wE(Z_j,\dots,Z_r)}}{\ab{\vfi_j\cdots\vfi_{r-1}}}\leq\sum_{j=1}^{r-5}\bigg(
\frac{10}{9}\bigg)^{r-j}(4u(\boldt))^3\bigg(\frac{1}{6}\bigg)^{r-j-5}\leq Cu^3(\boldt).\]
An application of Lemma \ref{Zj} yields
\[\frac{1}{\vfi_r}=\frac{1}{1-(1-\vfi_r)}=1+(1-\vfi_r)+\theta C\ab{1-\vfi_r}^2=-\EE Z_r+\theta C(k)u(\boldt).\]
The first and the third part of the  lemma follows by routinely applying Lemma \ref{Zj} and therefore we omit the detailed proof. For example,
\begin{eqnarray*}\EE \ab{Z_rZ_{r-1}Z_{r-2}Z_{r-3}}&\leq& 2\EE\ab{Z_1Z_2Z_3}\\
&=&16\sum_{l,j,m=1}^kp_{ljm}\ab{\sin(t_l/2)\sin(t_j/2)\sin(t_m/2)}\leq 16r_3(\boldt),\end{eqnarray*}
%\[\ab{\wE(Z_r,Z_{r-1})}\ab{1/\vfi_{r-1}-1}\leq  C(k)\sum_{j,m=1}^k\ab{p_{jm}-p_jp_m}\ab{\sin(t_j/2)\sin(t_m/2)}
%\sum_{l=1}^kp_l\ab{\sin(t_l/2)}\leq C(k)r_3(\boldt),\]
etc. For the proof of second and fourth expansions, observe that since all coordinates of $X_j$ are bounded by unity, we have
\[\ab{Z_j}\leq\ab{(\boldt,X_j)}\leq \sum_{j=1}^k\ab{t_j}.\]
Hence,
\[\EE \ab{Z_rZ_{r-1}Z_{r-2}Z_{r-3}}\leq \sum_{j=1}^k\ab{t_j}\EE\ab{Z_1Z_2Z_3}\leq \frac{1}{8}r_4(\boldt).\]
Here we have used also inequality $\ab{\sin(t_j/2)}\leq \ab{t_j}/2$.
We also apply the trivial expansion
\[(\ee^{\ii t_j}-1)(\ee^{\ii t_m}-1)(\ee^{\ii t_l}-1)=\ii^3 t_jt_mt_l+2\theta (\ab{t_l}+\ab{t_j}+\ab_{t_m})=\ii^3t_jt_mt_l+2\theta\sum_{l=1}^k\ab{t_l},\]
so that
\[\wE(Z_r,Z_{r-1},Z_{r-2})=\ii^3\sum_{l,j,m=1}^k(p_{ljm}-2p_{lj}p_m+p_lp_jp_m)+\theta C(k)r_4(\boldt).\]

For $k=2,3,4$ the proof is similar, the only difference being finite number of estimated summands in Heinrich's expansion. \qed

Let next
\[g_1(\boldt):= \Exponent{\EE Z_1-(\EE Z_1)^2/2},\quad
g_m(\boldt):=\Exponent{\EE Z_m-(\EE Z_m)^2/2+\wE(Z_m,Z_{m-1})}.\]

\begin{lemma}\label{gm} Assume the conditions in (\ref{a1a2}) hold. Then, for  $1 \leq r \leq n$,
\begin{eqnarray*}
g_r(\boldt)&=&1+\EE Z_r+\wE(Z_r,Z_{r-1})\\
&&+\theta C(k)\sum_{l,j,m=1}^k(p_{lj}p_m+p_lp_jp_m)\ab{\sin(t_l/2)\sin{t_j/2}\sin(t_m/2)},\\
g_r(\boldt)&=&1+\EE Z_r+\wE(Z_r,Z_{r-1})+\frac{\ii^3}{3}\sum_{l,j,m=1}^k(3p_{lj}p_m-4p_lp_jp_m)t_lt_jt_m\\
&&+\theta
C(k)\sum_{i=1}^k\sum_{l,j,m=1}^k(p_{lj}p_m+p_lp_jp_m)\ab{t_lt_jt_m},\\
\ab{g_r(\boldt)}&\leq& \exponent{-u(\boldt)}.
\end{eqnarray*}
\end{lemma}

\emph{Proof.} Observe that, by Lemma \ref{Zj} (see also the proof of Lemma \ref{C1}):
\begin{eqnarray*}\ab{g_r(\boldt)}&\leq&\exponent{-2u(\boldt)+\ab{\wE(Z_r,Z_{r-1})}+0.5\ab{\EE Z_r}^2}\\
&\leq&
\exponent{-2u(\boldt)+0.428 u(\boldt)+u(\boldt)/72}\leq\exponent{-u(\boldt)}.\end{eqnarray*}
Note also that the same estimate holds for $r=1$.
The Taylor series expansion gives us
\begin{eqnarray*}
g_r(\boldt)&=&1+\EE Z_r+\wE(Z_r,Z{r-1})(1+\EE Z)-\frac{1}{3}\EE^3 Z_r\\
&&+\theta C(\ab{\wE(Z_r,Z_{r-1})}\ab{\EE Z_r}
+\ab{\wE(Z_r,Z_{r-1})}^2+\ab{\EE Z_r}^4).
\end{eqnarray*}
The rest of the arguments is similar  to the proof of the previous lemma and, therefore, omitted.\qed
%PROOFA

%%%%%%%%%%%%%%%%%%%%%%%%%%%%%%%%%%%%%%%%%%%%%%%%%%%%%%%%%%%%%%%%%%%%%%%%%%%%%%%%%%%%%%%%%%%%%%%%%%%%%%%%%
%%%%%%%%%%%%%%%%%%%%%%%%%%%%%%%%%%%%%%%%%%%%%%%%%%%%%%%%%%%%%%%%%%%%%%%%%%%%%%%%%%%%%%%%%%%%%%%%%%%%%%%%%
%%%%%%%%%%%%%%%%%%%%%%%%%%%%%%%%%%%%%%%%%%%%%%%%%%%%%%%%%%%%%%%%%%%%%%%%%%%%%%%%%%%%%%%%%%%%%%%%%%%%%%%%%%%
%%%%%%%%%%%%%%%%%%%%%%%%%%%%%%%%%%%%%%%%%%%%%%%%%%%%%%%%%%%%%%%%%%%%%%%%%%%%%%%%%%%%%%%%%%%%%%%%%%%%%%%%%%%%5
%%%%%%%%%%%%%%%%%%%%%%%%%%%%%%%%%%%%%%%%%%%%%%%%%%%%%%%%%%%%%%%%%%%%%%%%%%%%%%%%%%%%%%%%%%%%%%%%%%%%%%%%%%%5
%%%%%%%%%%%%%%%%%%%%%%%%%%%%%%%%%%%%%%%%%%%%%%%%%%%%%%%%%%%%%%%%%%%%%%%%%%%%%%%%%%%%%%%%%%%%%%%%%%%%%%%%%%%%5

\section{Proofs}

\emph{Proofs of Theorems \ref{P} and \ref{G}.} Let
\[\pi_j(\boldt):=\exponent{\EE Z_j}=1+\EE Z_j +\theta C(k)\ab{\EE Z_j}^2.\]
 For simplicity, write $\pi_j= \pi_j(t)$. Using Lemmas \ref{Zj}, \ref{expofi} and \ref{vfiskl}, we have
\begin{eqnarray*}
\ab{\w F_n(\boldt)-\w\Pois (\bm{\lambda})}&=&\Ab{\prod_{j=1}^n\vfi_j-\prod_{j=1}^n\pi_j}\leq
\sum_{j=1}^n\ab{\vfi_j-\pi_j}\prod_{m=1}^j\ab{\vfi_j}\prod_{m=j+1}^n\ab{\pi_j}\\
&\leq& C(k)(r_1(\boldt)+\ab{\EE Z_j}^2)
\ee^{-n0.05u(\boldt)}
\leq C(k)r_1(\boldt)\ee^{-0.05nu(\boldt)}.
\end{eqnarray*}
 Here we have used also the trivial estimate $\exponent{0.04u(\boldt)}\leq C(k)$.
 Consequently applying Lemma \ref{locinv} and a lemma of \cite{Sho77}, we complete the proof for local and $\ell_2$ norms. Next observe that for $\alpha\geq 2$, the following trivial estimate holds:
 \begin{equation}
 \norm{M}_\alpha\leq \big(\norm{M}_\infty\big)^{(\alpha-2)/\alpha}\big(\norm{M}_2\big)^{2/\alpha}.
 \label{lpl2}
  \end{equation}

 The proof of Theorem \ref{G} is very similar. One needs
 to check that from Lemmas \ref{vfiskl} and \ref{gm} it follows that  $\ab{\vfi_m(\boldt)-g_m(\boldt)}\leq r_3(\boldt)$, for all $m=2,\dots,n$. It can be directly verified that the same estimate holds also for $m=1$. For the estimate of $\ab{g_{j+1}(\boldt)\cdots g_n(\boldt)}$ one should apply the last estimate of Lemma \ref{gm}. \hspace{.1cm} \hbox{\qed}

%%%%%%%%%%%%%%%%%%%%%%%%%%%%%%%%%%%%%%%%%%%%%%%%%%%%%%%%%%%%%%%%%%%%%%%%%%%%%%%%%%%%%%%%%%%%%%%%%%%%%%%%%%%%%%%%%%%%%%%%
%%%%%%%%%%%%%%%%%%%%%%%%%%%%%%%%%%%%%%%%%%%%%%%%%%%%%%%%%%%%%%%%%%%%%%%%%%%%%%%%%%%%%%%%%%%%%%%%%%%%%%%%%%%%%%%%%%%%%%%%%%
\emph{Proof of Theorem \ref{P1}.} It
suffices to estimate the difference of corresponding Fourier transforms.
 We again apply Lemmas \ref{Zj}, \ref{vfiskl}. By Bergstr\"om identity (see \cite{Ce16}, p. 17),
\begin{eqnarray*}
\lefteqn{\Ab{\prod_{m=1}^n\vfi_m-\prod_{m=1}^n\pi_m-\sum_{m=1}^n(\vfi_m-\pi_m)\prod_{j\ne m}^n\pi_j}}\hskip 1cm\\
&\leq&
\sum_{r=1}^n\ab{\vfi_r-\pi_r}\sum_{m=1}^{r-1}\ab{\vfi_m-\pi_m}\prod_{j=r+1}^n\ab{\vfi_j}\prod_{j\ne m}^n\ab{\pi_j}\\
&\leq& C(k) \exponent{-0.05nu(\boldt)}(nr_1(\boldt))^2.
\end{eqnarray*}
Next observe that $\ab{\pi_j-1}\leq C(k)\EE\ab{Z_m}\leq C(k)\sum_{r=1}^np_r\ab{\sin(t_r/2)}$. Therefore,
\[\Ab{\sum_{m=1}^n(\vfi_m-\pi_m)\prod_{j\ne m}^n\pi_j(1-\pi_m)}\leq C(k)\ee^{-0.05nu(\boldt)} r_1(\boldt)n\sum_{r=1}^np_r\ab{\sin(t_r/2)}.\]
It is easy to check that
\[\pi_m=1+\EE Z_m+0.5(\EE Z_k)^2+\theta C(k)\EE\ab{Z_k}^3.\]
Applying Lemma \ref{vfiskl}, we get
\[\prod_{j=1}^n\ab{\pi_j}\sum_{m=1}^n\ab{\vfi_m-\pi_m+0.5(\EE Z_m)^2-\wE(Z_m,Z_{m-1})}\leq C(k)\ee^{-0.05nu(\boldt)}r_3(\boldt).\]
Collecting all the estimates given above and applying Lemmas \ref{locinv},   a lemma of \cite{Sho77} and (\ref{lpl2}),  we complete the proof of theorem. \qed

%%%%%%%%%%%%%%%%%%%%%%%%%%%%%%%%%%%%%%%%%%%%%%%%%%%%%%%%%%%%%%%%%%%%%%%%%%%%%%%%%%%%%%%%%%%%%%%%%%%%%%%%%%%%%%%%%%%%%%%%
%%%%%%%%%%%%%%%%%%%%%%%%%%%%%%%%%%%%%%%%%%%%%%%%%%%%%%%%%%%%%%%%%%%%%%%%%%%%%%%%%%%%%%%%%%%%%%%%%%%%%%%%%%%%%%%%%%%%%%%%%%

\emph{Proof of Theorem \ref{LB}.} First we deal with approximation by signed compound Poisson measure $G$. We apply Lemma \ref{lblem} with
$\bolda=(n-1)(\boldt,\boldp):=(n-1)(p_1t_1+p_2t_2+\cdots+p_kt_k)$ and $\beta_j=b\sqrt{np_j}$, $j=1,2,\dots,n$.
From Lemma \ref{Zj} it follows that
\begin{eqnarray*}
\ab{\vfi_m(\boldt)\ee^{-\ii(\boldt,\boldp)}-1}&\leq&\ab{\vfi_m(\boldt)-\pi_m(\boldt)}\ab{\ee^{-\ii(\boldt,\boldp)}}\\
&&
+\Ab{\Exponent{\sum_{r=1}^kp_r(\ee^{\ii t_r}-1-\ii t_r)}-1}\\
&\leq& C(k)r_1(\boldt)+C(k)\sum_{r=1}^kp_rt_r^2\leq C(k)\sum_{r=1}^kp_rt_r^2.
\end{eqnarray*}
Here we for the last step we argued as in (\ref{vai}).  Similarly, we establish that the same estimate holds for
$\ab{g_m(\boldt)\ee^{-\ii(\boldt,\boldp)}-1}$.
Therefore, for any $m=1,2,\dots,n$,
  \[\Ab{\prod_{j=1}^{m-1}\vfi_j(\boldt_{\boldbeta})\prod_{j=m+1}^n g_j(\boldt_{\boldbeta})\ee^{-\ii(\bolda,\boldt_{\boldbeta})}-1}
  \leq C(k)n\sum_{r=1}^kp_r\frac{t^2_r}{b^2np_r}=\frac{C(k)}{b^2}\sum_{r=1}^kt^2_r.\]
 Consequently applying Lemmas \ref{vfiskl} and \ref{gm} and using the above estimate, we get
\begin{eqnarray}
 \lefteqn{\Ab{\ee^{-\ii(\boldt_{\boldbeta},\boldp)}\bigg(\prod_{j=1}^n\vfi_j(\boldt_{\boldbeta})-\prod_{j=1}^n g_j(\boldt_{\boldbeta})\bigg)-
 \sum_{j=1}^n(\vfi_j(\boldt_{\boldbeta})-g_j(\boldt_{\boldbeta}))}}\hskip 1cm\nonumber\\
 &\leq& \sum_{m=1}^n
 \ab{\vfi_j(\boldt_{\boldbeta})-g_j(\boldt_{\boldbeta})}{C(k)}{b^2} \sum_{r=1}^kt_r^2
 \leq \frac{C(k)}{b^2} \sum_{r=1}^k\ab{t_r} r_4(\boldt_{\boldbeta}).
 \label{zem1}
 \end{eqnarray}
Expanding $g_1(\boldt)$ in Taylor series and noting that the remainder term is smaller than $r_4(\boldt)$, we obtain
\[\vfi_1(\boldt)-g_1(\boldt)=\ii^3\sum_{r,j,m=1}^kp_rp_jp_mt_rt_jt_m+\theta C(k)r_4(\boldt).\]
Similarly
\[ \vfi_2(\boldt)-g_2(\boldt)=\frac{\ii^3}{3}\sum_{r,j,m=1}^k(6p_{rj}p_m-5p_rp_jp_m)t_rt_jt_m+\theta C(k)r_4(\boldt)\]
and from Lemmas \ref{vfiskl} and \ref{gm}, for $i=3,4,\dots,n$,
\[\vfi_i(\boldt)-g_i(\boldt)=\frac{\ii^3}{3}\sum_{r,j,m=1}^k(3p_{rjm}-12p_{rj}p_m+10p_rp_jp_m)t_rt_jt_m+\theta C(k)r_4(\boldt).\]
Therefore,
\begin{eqnarray}
\lefteqn{\sum_{i=1}^n(\vfi_i(\boldt_{\boldbeta})-g_i(\boldt_{\boldbeta}))=\frac{\ii^3}{3b^3\sqrt{n}}\sum_{r,j,m=1}^k
\frac{3p_{rjm}-12p_{rj}p_m+10p_rp_jp_m}{\sqrt{p_rp_jp_m}}t_rt_jt_m}\hskip 1cm\nonumber\\
&&+\theta \frac{1}{3b^3n\sqrt{n}}\Ab{\sum_{r,j,m=1}^k\frac{-6p_{rjm}+30p_{rj}p_m-22p_rp_jp_m}{\sqrt{p_rp_jp_m}}t_rt_jt_m}\nonumber\\
&&+\theta C(k)n
r_4(\boldt_{\boldbeta})
. \label{zem2}
\end{eqnarray}

Combining the estimates in (\ref{zem1}) and (\ref{zem2}) and observing that integral of $\ab{t_l}^s\psi_l(t_l)$ is bounded by absolute constant, we
can write
\begin{equation}
\ab{V(\bolda,\boldbeta)}\geq \frac{1}{3b^3\sqrt{n}}\bigg\vert\int_{-\infty}^\infty\cdots\int_{-\infty}^\infty  \prod_{l=1}^k\psi_l(t_l)
\sum_{r,j,m=1}^k v_{rjm}t_rt_jt_m\bigg\vert
-
\frac{C(k)\veps_6}{b^3\sqrt{n}}.
\label{zem3}
\end{equation}
Here
\[v_{rjm}:=\frac{3p_{rjm}-12p_{rj}p_m+10p_rp_jp_m}{\sqrt{p_rp_jp_m}}.\]

  Now comes the tricky part, since all integrals with odd powers of $t_l$ are equal zero. We can write
  \begin{eqnarray*}
  \sum_{r,j,m=1}^k v_{rjm}&=&\sum_{r=1}^kt_r^3v_{rrr}+\sum_{r\ne m}^k t_r^2t_m(v_{rrm}+b_{rmr}+v_{mrr})\\
  &&+
  \sum_{r\ne m\ne j}^k t_rt_j t_m(v_{rmj}+b_{rjm}+v_{mrj}+v_{mrl}+v_{jmr}+v_{jrm}).
  \end{eqnarray*}
Let us choose $\psi_m(t_m)=t_m\ee^{-t^2_m/2}$ and  $\psi_t(t_r)=\ee^{-t_r^2/2}$, for all $r\ne m$.
Then, after integration, absolute value in (\ref{zem3}) is equal to
\[(\sqrt{2\pi})^k\Ab{3v_{mmm}+\sum_{r=1,r\ne m}^k(v_{rmm}+v_{mrm}+v_{mmr})}.\]
By taking different $m$ we obtain $k$ such integrals. Now, let us assume that $\psi_m=t_m\ee^{-t^2_m/2}$, $\psi_r=t_r\ee^{-t^2_r/2}$,
$\psi_j=t_j\ee^{-t^2_j/2}$ and all other $\psi_l(t_l)=\ee^{-t_l^2/2}$.
Then, after integration, absolute value in (\ref{zem3}) is equal to
\[(\sqrt{2\pi})^k\ab{v_{rmj}+v_{rjm}+v_{mrj}+v_{mrl}+v_{jmr}+v_{jrm}}.\]
 After taking all possible different combinations $r,m,j$, we arrive at the fact that absolute value  in (\ref{zem3}) can be taken equal to maximum of all these $N=k+k(k-1)(k-2)/6$ estimates. Next observe that for any numbers $B_1,\dots,B_N$, we have
  \[\max_{1\leq j\leq N}\ab{B_j}\geq \frac{1}{N}\sum_{j=1}^N\ab{B_j}\geq \frac{1}{N}\Ab{\sum_{j=1}^N B_j}.\]
Therefore,
\[\bigg\vert\int_{-\infty}^\infty\cdots\int_{-\infty}^\infty  \prod_{m=1}^k\psi_m(t_m)
\sum_{r,j,m=1}^kv_{rjm}t_rt_jt_m\bigg\vert\geq \frac{(\sqrt{2\pi})^k}{N}\bigg\vert \sum_{r,j,m=1}^kv_{rjm}
+2 \sum_{m=1}^k v_{mmm}\bigg\vert.\]
Collecting all the  relevant estimates, we complete the proof for approximation $G$. The estimates for Poisson approximation are obtained by the similar arguments. Note  that as we need to integrate sums of the form $\sum_{r,m=1}^k w_{rm}t_rt_m$, the choice of $\psi_r(t_r)$, $\psi_m(t_m)$ allows to estimate corresponding integral by
\begin{eqnarray*}
\lefteqn{(\sqrt{2\pi})^k\max\bigg\{\Ab{\sum_{m=1}^k w_{mm}},\ab{w_{12}},\ab{w_{13}},\dots,\ab{w_{k-1,k}}\bigg\}}
\hskip 2cm\\
&\geq& (\sqrt{2\pi})^k\frac{2}{k(k+1)}\Ab{\sum_{r,m=1}^k w_{rm}}.
 \hspace{2.5cm} \hbox{\qed}
\end{eqnarray*}

\emph{Proof of Proposition \ref{prop}}. By triangle inequality,
\[\norm{F_n-\Pois(\bm{\lambda})}_\infty\leq \norm{F_n-\Pois(\bm{\lambda})-A_1}_\infty+\norm{A_1}_\infty.\]
Observe that
\begin{eqnarray*}
\lefteqn{-\frac{n}{2}\bigg(\sum_{j=1}^2p_j(\ee^{\ii t_j}-1)\bigg)^2+
(n-1)\sum_{j,m=1}^2(p_{jm}-p_jp_m)(\ee^{\ii t_j}-1)(\ee^{\ii t_m}-1)}\hskip 1cm\\
&=&n\delta^2q^3(\ee^{\ii t_1}-1)^2(1+o(1)) +n(1-\delta)^2\bar{q}^3(\ee^{\ii t_2}-1)(1+o(1))\\
&&-n\delta(1-\delta)q^2\bar{q}^2(\ee^{\ii t_1}-1)(\ee^{\ii t_2}-1).
\end{eqnarray*}
It remains to apply Lemmas \ref{locinv}, \ref{shorgin} and the obvious inequality
\begin{align*}
\sqrt{\delta(1-\delta)}q\bar{q}\leq\frac{\delta q^2+(1-\delta)\bar{q}^2}{2} \leq
\frac{\delta q+(1-\delta)\bar{q}}{2}. \hspace{2cm} \hbox{\qed}
\end{align*}
%%%%%%%%%%%%%%%%%%%%%%%%%%%%%%%%%%%%%%%%%%%%%%%%%%%%%%%%%%%%%%%%%%%%%%%%%%%%%%%%%%%%%%%%%%%%%%%%%%%%%%%%%%%%%%%%%%%%%%%
%

%\section{Section title}
%\label{sec:1}
%Text with citations \cite{RefB} and \cite{RefJ}.
%\subsection{Subsection title}
%\label{sec:2}
%as required. Don't forget to give each section
%and subsection a unique label (see Sect.~\ref{sec:1}).
%\paragraph{Paragraph headings} Use paragraph headings as needed.
%\begin{equation}
%\end{equation}

%\begin{figure}
%  \includegraphics{example.eps}
%\end{figure}
%

%\begin{table}
%\begin{tabular}{lll}
%\hline\noalign{\smallskip}
%\end{tabular}
%\end{table}

\subsection*{Acknowledgements}
This paper was completed during the first author's stay (January-February, 2020) at the
Department of Mathematics, Indian Institute of Technology Bombay. The first author
 would like to thank  the Department and the Institute for the
 invitation and the hospitality. We are
grateful to the referee for useful remarks, which lead to the improvement of the paper.




\end{document}